\documentclass[11pt,draft]{article}

\usepackage{a4,amsmath,amssymb,amsthm,cite}

\DeclareMathOperator{\var}{var}

\newtheorem{theorem}{Theorem}

\newtheorem{lemma}[theorem]{Lemma}

\newtheorem{remark}[theorem]{Remark}

\begin{document}

\title{\textbf{Chains and anti-chains}\\
\textbf{in the lattice of epigroup varieties}\thanks{The work was partially
supported by the Russian Foundation for Basic Research (grant
No.\,09-01-12142) and the Federal Education Agency of the Russian Federation
(project No.~2.1.1/3537).}}

\author{D.\,V.\,Skokov and B.\,M.\,Vernikov}

\date{}

\maketitle

\begin{abstract}
Let $\mathcal E_n$ be the variety of all epigroups of index $\le n$. We prove
that, for an arbitrary natural number $n$, the interval $[\mathcal E_n,
\mathcal E_{n+1}]$ of the lattice of epigroup varieties contains a chain
isomorphic to the chain of real numbers with the usual order and an
anti-chain of the cardinality continuum.

\emph{Key words}: epigroup, variety, lattice of subvarieties.

\emph{AMS Subject Classification}: primary 20M07, secondary 08B15.
\end{abstract}

A semigroup $S$ is called an \emph{epigroup} if for any element $x$ of $S$
some power of $x$ lies in some subgroup of $S$. For an element $a$ of a given
epigroup, let $e_a$ be the unit element of the maximal subgroup $G$ that
contains some power of $a$. It is known that $ae_a=e_aa$ and this element
lies in $G$. We denote by $\overline a$ the element inverse to $ae_a$ in $G$.
This element is called the \emph{pseudo-inverse} of $a$. The mapping $a
\longmapsto\overline a$ defines a unary operation on an epigroup. The idea to
treat epigroups as unary semigroups (that is semigroups with an additional
unary operation of pseudo-inversion) was promoted by Shevrin in \cite
{Shevrin-94}. A systematic overview of the material accumulated in the theory
of epigroups by the beginning of the 2000s was given in the survey \cite
{Shevrin-05}.

By epigroup variety we mean a variety of epigroups treated just as unary
semigroups. Results about epigroup varieties that are known so far mainly
concern with equational and structural aspects (see corresponding results in
\cite{Shevrin-94,Shevrin-05}). As to considerations of the varietal lattices,
there are only a few results about such a type (see Sections~2 and~3 in the
recent survey \cite{Shevrin-Vernikov-Volkov-09}). In \cite{Shevrin-94}
several open questions about lattices of epigroup varieties were formulated;
some of them are reproduced in \cite{Shevrin-05} and \cite
{Shevrin-Vernikov-Volkov-09}. The aim of this note is to answer one of these
questions and obtain an information closely related with one more of them.

An epigroup $S$ has \emph{index} $n$ if the $n$th power of every element of
$S$ lies in some of its subgroups and $n$ is the least number with this
property. The class of all epigroups of index $\le n$ is denoted by $\mathcal
E_n$. For each $n$, the class $\mathcal E_n$ is known to be a variety of
epigroups; it is given by the identities
$$(xy)z=x(yz),\ x\overline x=\overline xx,\ x\overline x\,^2=\overline x,\
x^{n+1}\overline x=x^n$$
(see \cite{Shevrin-94}). The chain $\mathcal E_1\subset\mathcal E_2\subset
\cdots\mathcal E_n\subset\cdots$ can be regarded as the ``spine'' of the
lattice of all epigroup varieties, since for any epigroup variety $\mathcal
V$ there exists $n$ such that $\mathcal{V\subseteq E}_n$.

The following questions have been formulated in \cite{Shevrin-94} and
repeated in \cite{Shevrin-05,Shevrin-Vernikov-Volkov-09}:

\smallskip

1) \emph{What are the order types of maximal chains in the intervals
$[\mathcal E_n,\mathcal E_{n+1}]$ of the lattice of epigroup varieties}?

\smallskip

2) \emph{What are the cardinalities of maximal anti-chains in these
intervals}?

\smallskip

The first question is still open. But the following theorem shows that the
intervals $[\mathcal E_n,\mathcal E_{n+1}]$ contain rather complicated
chains.

\begin{theorem}
\label{chain}
For an arbitrary natural number $n$, the interval $[\mathcal E_n,\mathcal
E_{n+1}]$ contains a chain isomorphic to the chain of real numbers with the
usual order.
\end{theorem}

Note that chains we construct in the proof of Theorem~\ref{chain} are not
maximal in the intervals of the kind $[\mathcal E_n,\mathcal E_{n+1}]$ (see
Remark \ref{not max} below).

The complete answer on the second question is given by the following

\begin{theorem}
\label{anti-chain}
For an arbitrary natural number $n$, the interval $[\mathcal E_n,\mathcal
E_{n+1}]$ contains an anti-chain of cardinality continuum.
\end{theorem}

There are two results that play the key role in the proof of both theorems.
The first of them was proved by Je\v{z}ek in \cite{Jezek-76}. To formulate
this result, we recall that a word $u$ is said to be \emph{applicable to a
word} $v$ if $v$ may be presented in the form $a\xi(u)b$ where $a$ and $b$
are (maybe empty) words, while $\xi$ is an endomorphism on the free semigroup
under a countably infinite alphabet. The mentioned result by Je\v{z}ek is
that there are a countably infinite set of semigroup words $\{w_i\mid i\in
I\}$ such that $w_i$ is not applicable to $w_j$ for any $i,j\in I$, $i\ne j$
and $x^2$ is not applicable to $w_i$ for any $i\in I$. For our aim, it is
convenient to enumerate these words by rational numbers. In what follows we
will refer to these words as to the words $Z_\alpha$ where $\alpha$ runs over
the set of all rational numbers. For each rational $\alpha$, the first letter
of $Z_\alpha$ will be denoted by $x_\alpha$.

To formulate the second result, we need some definitions and notation. A pair
of identities $wx=xw=w$ where the letter $x$ does not occur in the word $w$
is usually written as the symbolic identity $w=0$. (This notation is
justified because a semigroup with the identities $wx=xw=w$ has a zero
element and all values of the word $w$ in this semigroup are equal to zero.)
An identity of the form $w=0$ as well as a variety given by identities of
such a form are called 0-\emph{reduced}.  A semigroup variety is called a
\emph{nil-variety} if it consists of nil-semigroups; this takes place if and
only if it satisfies the identity $x^n=0$ for some $n$. It is evident that
every 0-reduced variety is a nil-variety. It is clear that every
nil-semigroup is an epigroup and every nil-variety of semigroups may be
considered as a variety of epigroups.

An element $x$ of a lattice $\langle L;\vee,\wedge\rangle$ is called \emph
{lower-modular} if
$$\forall\,y,z\in L\colon\quad x\le y\longrightarrow(z\vee x)\wedge y=(z
\wedge y)\vee x\ldotp$$
\emph{Upper-modular} elements are defined dually. It was verified in \cite
[Corollary~3]{Vernikov-Volkov-88} that a 0-reduced semigroup variety is a
lower-modular element of the lattice of all semigroup varieties. The proof of
this fact given in \cite{Vernikov-Volkov-88} is based on the following two
ingredients: 1)~the fully invariant congruence on the free semigroup
corresponding to a 0-reduced variety has exactly one non-singleton class;
2)~an equivalence relation $\pi$ on a set $S$ has at most one non-singleton
class if and only if $\pi$ is an upper-modular element of the equivalence
lattice of $S$ (this observation was checked in \cite[Proposition~3]
{Vernikov-Volkov-88}). It is evident that these arguments are applicable for
epigroup varieties as well. Thus we have

\begin{lemma}
\label{0-red is lmod}
A \textup0-reduced epigroup variety is a lower-modular element of the lattice
of all epigroup varieties.\qed
\end{lemma}

A semigroup variety given by an identity system $\Sigma$ is denoted by $\var
\Sigma$.

Now we are ready to prove both theorems.

\emph{Proof of Theorem} \ref{chain}. Let $n$ be a natural number and $\xi$ a
real number. Put
$$\mathcal C_\xi^n=\var\,\{x^{n+1}=x_\alpha^{n-1}Z_{\alpha}=0\mid \alpha\ge
\xi\}$$
(if $n=1$ then $x_\alpha^0$ is the empty word) and $\mathcal D_\xi^n=\mathcal
E_n\vee\mathcal C_\xi^n$. It is clear that $\mathcal C_\xi^n\subseteq\mathcal
E_{n+1}$, whence $\mathcal D_\xi^n\in[\mathcal E_n,\mathcal E_{n+1}]$. Let
now $\xi_1$ and $\xi_2$ be real numbers with $\xi_1\le\xi_2$. Then $\mathcal
C_{\xi_1}^n\subseteq\mathcal C_{\xi_2}^n$ and therefore $\mathcal D_{\xi_1}^n
\subseteq\mathcal D_{\xi_2}^n$. To prove Theorem~\ref{chain}, it suffices to
verify that $\mathcal D_{\xi_1}^n\ne\mathcal D_{\xi_2}^n$ whenever $\xi_1\ne
\xi_2$. Arguing by contradiction, suppose that $\xi_1<\xi_2$ (and therefore
$\mathcal C_{\xi_1}^n\subset\mathcal C_{\xi_2}^n$) but $\mathcal D_{\xi_1}^n=
\mathcal D_{\xi_2}^n$ (see Fig.~\ref{fig chain}).

Note that all varieties of the kind $\mathcal C_\xi^n$ are 0-reduced.
Further, for any $\xi$, the variety $\mathcal E_n\wedge\mathcal C_\xi^n$ is a
nil-variety of index $\le n$, whence it satisfies the identity $x^n=0$.
Therefore
\begin{equation}
\label{eq chain}
\mathcal E_n\wedge\mathcal C_{\xi_2}^n\subseteq\mathcal C_{\xi_1}^n\ldotp
\end{equation}
We have
\begin{align*}
\mathcal C_{\xi_1}^n&=(\mathcal E_n\wedge\mathcal C_{\xi_2}^n)\vee\mathcal
C_{\xi_1}^n&&\text{by \eqref{eq chain}}\\
&=(\mathcal E_n\vee\mathcal C_{\xi_1}^n)\wedge\mathcal C_{\xi_2}^n&&\text{by
Lemma \ref{0-red is lmod}}\\
&=\mathcal D_{\xi_1}^n\wedge\mathcal C_{\xi_2}^n&&\text{by the definition
of}\ \mathcal D_{\xi_1}^n\\
&=\mathcal D_{\xi_2}^n\wedge\mathcal C_{\xi_2}^n&&\text{because}\ \mathcal
D_{\xi_1}^n=\mathcal D_{\xi_2}^n\\
&=\mathcal C_{\xi_2}^n&&\text{by the definition of}\ \mathcal D_{\xi_2}^n
\ldotp
\end{align*}
Thus $\mathcal C_{\xi_1}^n=\mathcal C_{\xi_2}^n$. A contradiction.\qed

Let $C=\{\mathcal D_\xi^n\mid\xi\in\mathbb R\}$. If $\xi\in\mathbb R$ then
$\mathcal E_n\ne\mathcal D_\xi^n$ because $\mathcal C_\xi^n\nsubseteq\mathcal
E_n$, and $\mathcal E_{n+1}\ne\mathcal D_\xi$ because $\mathcal{D_\xi\subset
D_\lambda\subseteq E}_{n+1}$ for any $\lambda\in\mathbb R$ with $\xi<
\lambda$. Thus, we may ajoin $\mathcal E_n$ [respectively $\mathcal E_{n+1}$]
as the least [the greatest] element to the chain $C$ and obtain a chain $C^*$
in $[\mathcal E_n,\mathcal E_{n+1}]$ with $C\subset C^*$. We have the
following

\begin{remark}
\label{not max}
The chain $C$ is not the maximal chain in the interval $[\mathcal E_n,
\mathcal E_{n+1}]$.\qed
\end{remark}

\emph{Proof of Theorem} \ref{anti-chain}. As in the proof of Theorem \ref
{chain}, let $n$ be a natural number and $\xi$ a real number. Now we put
$$\mathcal A_\xi^n=\var\,\{x^{n+1}=x_\alpha^{n-1}Z_{\alpha}=0\mid\xi-1<
\alpha<\xi+1\}$$
and $\mathcal B_{\xi}^n=\mathcal E_n \vee\mathcal A_\xi^n$. It is clear that
$\mathcal A_\xi^n\subseteq\mathcal E_{n+1}$ and $\mathcal B_\xi^n\in[\mathcal
E_n,\mathcal E_{n+1}]$. Let $\xi_1$ and $\xi_2$ be different real numbers.
Then the varieties $\mathcal A_{\xi_1}^n$ and $\mathcal A_{\xi_2}^n$ are
non-comparable. To prove Theorem~\ref{anti-chain}, it suffices to verify that
the varieties $\mathcal B_{\xi_1}^n$ and $\mathcal B_{\xi_2}^n$ are
non-comparable too. Arguing by contradiction, suppose that $\mathcal
B_{\xi_2}^n\subseteq\mathcal B_{\xi_1}^n$ (see Fig.~\ref{fig anti-chain}).

\begin{table}[tbh]
\begin{center}
\begin{tabular}{c@{\hskip75pt}c}
\unitlength=1mm
\linethickness{0.4pt}
\begin{picture}(23,33)
\put(1,10){\line(1,-1){10}}
\put(1,10){\line(1,2){10}}
\put(11,0){\line(1,1){10}}
\put(11,30){\line(1,-1){10}}
\put(21,10){\line(0,1){10}}
\put(1,10){\circle*{1.33}}
\put(11,0){\circle*{1.33}}
\put(11,30){\circle*{1.33}}
\put(21,10){\circle*{1.33}}
\put(21,20){\circle*{1.33}}
\put(23,10){\makebox(0,0)[lc]{$\mathcal C_{\xi_1}^n$}}
\put(23,20){\makebox(0,0)[lc]{$\mathcal C_{\xi_2}^n$}}
\put(11,33){\makebox(0,0)[cc]{$\mathcal D_{\xi_1}^n=\mathcal D_{\xi_2}^n$}}
\put(0,10){\makebox(0,0)[rc]{$\mathcal E_n$}}
\end{picture}
\refstepcounter{figure}
\label{fig chain}&
\unitlength=1mm
\linethickness{0.4pt}
\begin{picture}(43,43)
\put(1,20){\line(1,-1){20}}
\put(1,20){\line(1,1){20}}
\put(11,30){\line(2,-1){20}}
\put(21,0){\line(1,1){20}}
\put(21,40){\line(2,-1){20}}
\put(31,10){\line(0,1){10}}
\put(31,20){\line(1,1){10}}
\put(41,20){\line(0,1){10}}
\put(1,20){\circle*{1.33}}
\put(11,30){\circle*{1.33}}
\put(21,0){\circle*{1.33}}
\put(21,40){\circle*{1.33}}
\put(31,10){\circle*{1.33}}
\put(31,20){\circle*{1.33}}
\put(41,20){\circle*{1.33}}
\put(41,30){\circle*{1.33}}
\put(43,20){\makebox(0,0)[lc]{$\mathcal A_{\xi_1}^n$}}
\put(43,30){\makebox(0,0)[lc]{$\mathcal A_{\xi_1}^n\vee\mathcal
A_{\xi_2}^n$}}
\put(39,8){\makebox(0,0)[cc]{$\mathcal A_{\xi_1}^n\wedge\mathcal
A_{\xi_2}^n$}}
\put(28,17){\makebox(0,0)[cc]{$\mathcal A_{\xi_2}^n$}}
\put(22,43){\makebox(0,0)[cc]{$\mathcal B_{\xi_1}^n$}}
\put(8,32){\makebox(0,0)[cc]{$\mathcal B_{\xi_2}^n$}}
\put(0,20){\makebox(0,0)[rc]{$\mathcal E_n$}}
\end{picture}
\refstepcounter{figure}
\label{fig anti-chain}\\
\small Figure \ref{fig chain}&\small Figure \ref{fig anti-chain}
\end{tabular}
\end{center}
\end{table}

Note that all varieties of the kind $\mathcal A_\xi^n$ are 0-reduced.
Further, the variety $\mathcal E_n\wedge(\mathcal A_{\xi_1}^n\vee\mathcal
A_{\xi_2}^n)$ is a nil-variety of index $\le n$, whence it satisfies the
identity $x^n=0$. Therefore,
\begin{equation}
\label{eq1 anti-chain}
\mathcal E_n\wedge(\mathcal A_{\xi_1}^n\vee\mathcal A_{\xi_2}^n)\subseteq
\mathcal A_{\xi_1}^n\ldotp
\end{equation}
Furthermore, $\mathcal B_{\xi_1}^n\supseteq\mathcal A_{\xi_1}^n$ and
$\mathcal B_{\xi_1}^n\supseteq\mathcal B_{\xi_2}^n\supseteq\mathcal
A_{\xi_2}^n$, whence
\begin{equation}
\label{eq2 anti-chain}
\mathcal B_{\xi_1}^n\supseteq\mathcal A_{\xi_1}^n\vee\mathcal A_{\xi_2}^n
\ldotp
\end{equation}
We have
\begin{align*}
\mathcal A_{\xi_1}^n&=\bigl(\mathcal E_n\wedge(\mathcal A_{\xi_1}^n\vee
\mathcal A_{\xi_2}^n)\bigr)\vee\mathcal A_{\xi_1}^n&&\text{by \eqref
{eq1 anti-chain}}\\
&=(\mathcal E_n\vee\mathcal A_{\xi_1}^n)\wedge(\mathcal A_{\xi_1}^n\vee
\mathcal A_{\xi_2}^n)&&\text{by Lemma \ref{0-red is lmod}}\\
&=\mathcal B_{\xi_1}^n\wedge(\mathcal A_{\xi_1}^n\vee\mathcal A_{\xi_2}^n)&&
\text{by the definition of}\ \mathcal B_{\xi_1}^n\\
&=\mathcal A_{\xi_1}^n\vee\mathcal A_{\xi_2}^n&&\text{by \eqref
{eq2 anti-chain}}\ldotp
\end{align*}
Thus $\mathcal A_{\xi_1}^n=\mathcal A_{\xi_1}^n\vee\mathcal  A_{\xi_2}^n$,
whence $\mathcal A_{\xi_2}^n\subseteq\mathcal A_{\xi_1}^n$. A contradiction.
\qed

\medskip

\textbf{Acknowledgements.} The authors would like to thank Professor
M.\,V.\,Vol\-kov for fruitful discussions.

\begin{flushleft}
Department of Mathematics and Mechanics, Ural State University, Lenina 51,

620083 Ekaterinburg, Russia

\smallskip

\emph{E-mail address}: \texttt{dskokov@yandex.ru, boris.vernikov@usu.ru}
\end{flushleft}

\end{document}